\documentclass[10pt]{article}
\usepackage{indentfirst, latexsym, amsmath,bm,amssymb}
\usepackage{tikz}
\usepackage{color,soul}
\usepackage{amsmath}
\usepackage{amssymb}
\usepackage{mathrsfs}
\usepackage{mathtools}
\usepackage{stmaryrd}
\usepackage{bbold}

\begin{document}

\newtheorem{theorem}{Theorem}[section]
\newtheorem{corollary}[theorem]{Corollary}
\newtheorem{definition}[theorem]{Definition}
\newtheorem{conjecture}[theorem]{Conjecture}
\newtheorem{question}[theorem]{Question}
\newtheorem{lemma}[theorem]{Lemma}
\newtheorem{remark}[theorem]{Remark}
\newtheorem{proposition}[theorem]{Proposition}
\newtheorem{example}[theorem]{Example}
\newenvironment{proof}{\noindent {\bf
Proof.}}{\rule{3mm}{3mm}\par\medskip}
\newcommand{\pp}{{\it p.}}
\newcommand{\de}{\em}

\allowdisplaybreaks[4]
\newtheorem {Problem} {Problem}[section]
\newtheorem {Theorem} [Problem]{Theorem}
\newtheorem {Lemma}[Problem]{Lemma}
\newtheorem{Conjecture}[Problem]{Conjecture}
\newtheorem {Corollary}[Problem]{Corollary}
\newenvironment {Proof}{\noindent {\bf Proof.}}{\hfill\ensuremath{\square}}
\newcommand*{\QEDB}{\hfill\ensuremath{\square}}

\newcommand{\JEC}{{\it Europ. J. Combinatorics},  }
\newcommand{\JCTB}{{\it J. Combin. Theory Ser. B.}, }
\newcommand{\JCT}{{\it J. Combin. Theory}, }
\newcommand{\JGT}{{\it J. Graph Theory}, }
\newcommand{\ComHung}{{\it Combinatorica}, }
\newcommand{\DM}{{\it Discrete Math.}, }
\newcommand{\ARS}{{\it Ars Combin.}, }
\newcommand{\SIAMDM}{{\it SIAM J. Discrete Math.}, }
\newcommand{\SIAMADM}{{\it SIAM J. Algebraic Discrete Methods}, }
\newcommand{\SIAMC}{{\it SIAM J. Comput.}, }
\newcommand{\ConAMS}{{\it Contemp. Math. AMS}, }
\newcommand{\TransAMS}{{\it Trans. Amer. Math. Soc.}, }
\newcommand{\AnDM}{{\it Ann. Discrete Math.}, }
\newcommand{\NBS}{{\it J. Res. Nat. Bur. Standards} {\rm B}, }
\newcommand{\ConNum}{{\it Congr. Numer.}, }
\newcommand{\CJM}{{\it Canad. J. Math.}, }
\newcommand{\JLMS}{{\it J. London Math. Soc.}, }
\newcommand{\PLMS}{{\it Proc. London Math. Soc.}, }
\newcommand{\PAMS}{{\it Proc. Amer. Math. Soc.}, }
\newcommand{\JCMCC}{{\it J. Combin. Math. Combin. Comput.}, }
\newcommand{\GC}{{\it Graphs Combin.}, }

\title{On the  spectral radius of  graphs  without a star forest\thanks{
This work is supported by  the National Natural Science Foundation of China (Nos. 11971311, 11531001), the Montenegrin-Chinese Science and Technology Cooperation Project (No.3-12)
}}

\author{Ming-Zhu Chen\thanks{E-mail: mzchen@hainanu.edu.cn}, A-Ming Liu,\thanks{E-mail: amliu@hainanu.edu.cn}\\
School of Science, Hainan University, Haikou 570228, P. R. China, \\
\and  Xiao-Dong Zhang$^{\dagger}$\thanks{Corresponding Author E-mail: xiaodong@sjtu.edu.cn}
\\School of Mathematical Sciences, MOE-LSC, SHL-MAC\\
Shanghai Jiao Tong University,
Shanghai 200240, P. R. China}

\date{}
\maketitle

\begin{abstract}
In this paper, we determine the  maximum spectral radius and all extremal graphs for  (bipartite) graphs  of order $n$ without  a star forest, extending  Theorem~1.4 (iii) and Theorem~1.5 for large $n$.  As a corollary, we determine the minimum least eigenvalue  of $A(G)$  and  all  extremal graphs for graphs of order $n$ without a star forest, extending  Corollary~1.6 for large $n$.
\end{abstract}

\section{Introduction}
 Let $G$ be an undirected simple graph with vertex set
$V(G)=\{v_1,\dots,v_n\}$ and edge set $E(G)$, where $n$ is called the order of $G$.
The \emph{adjacency matrix}
$A(G)$ of $G$  is the $n\times n$ matrix $(a_{ij})$, where
$a_{ij}=1$ if $v_i$ is adjacent to $v_j$, and $0$ otherwise. The  \emph{spectral radius} of $G$ is the largest eigenvalue of $A(G)$, denoted by $\rho(G)$.  The least eigenvalue of $A(G)$ is denoted by $\rho_n(G)$.
For $v\in V(G)$,  the \emph{neighborhood} $N_G(v)$ of $v$  is $\{u: uv\in E(G)\}$ and the \emph{degree} $d_G(v)$ of $v$  is $|N_G(v)|$.
We write $N(v)$ and $d(v)$ for $N_G(v)$ and $d_G(v)$ respectively if there is no ambiguity.
 Denote by $\Delta(G)$ the maximum degree of $G$.  Let $S_{n-1}$ be a star of order $n$. The \emph{center} of a star is the vertex of maximum degree in the star.   The \emph{ centers } of a star forest are the centers of the stars in the star forest. A graph $G$ is \emph{$H$-free} if   it  does not contain $H$ as a subgraph.
 For two vertex disjoint graphs $G$ and $H$,  we denote by  $G\cup H$ and  $G\nabla H$  the \emph{union} of $G$ and $H$,
and the \emph{join} of $G$ and $H$ which is obtained by joining every vertex of $G$ to every vertex of $H$, respectively.
Denote by $kG$  the  the union of $k$ disjoint copies of $G$.
For graph notation and terminology undefined here,  readers are referred to \cite{BM}.

Recall that the problem of maximizing the number of edges over all graphs without  fixed subgraphs is one of the cornerstones of graph theory.

\begin{Problem}\label{P0}
Given a graph $H$, what is the maximum number of edges of a graph $G$ of order $n$ without $H$ $?$
\end{Problem}

Many instances of Problems~\ref{P0} have been solved. For example,
 Lidick\'{y}, Liu, and Palmer \cite{LLP}  determined the maximum number of edges of  graphs without a  forest forest  if the order of a  graph is  sufficiently large.

\begin{Theorem}\label{thm-edge}\cite{LLP}
Let $F=\cup_{i=1}^k S_{d_i}$ be a star forest with $k\geq2$ and $d_1\geq\cdots \geq d_k\geq2$.
If $G$ is an $F$-free graph of  sufficiently large order  $n$, then
 $$e(G)\leq \max_{1\leq i\leq k}\bigg\{(i-1)(n-i+1)+\binom{i-1}{2}+\bigg\lfloor\frac{(d_i-1)(n-i+1)}{2}\bigg\rfloor\bigg\}.$$
\end{Theorem}

In spectral extremal graph theory, a similar central problem  is
of the following type:
\begin{Problem}\label{P1}
Given a graph $H$, what is the maximum  $\rho(G)$ of a graph $G$ of order $n$ without $H$$?$
\end{Problem}

Many instances of Problem~\ref{P1} have been solved, for example, see \cite{Chen2019-2,Cioaba2019,Gao2019,Nikiforov2010,Tait2019,Tait2017,Zhai2020}. In addition,  if $H$ is  a linear forest, Problem~\ref{P1} was solved in \cite{Chen2019}. For $H=kP_3$, the bipartite version  of Problem~\ref{P1} was also proved in  \cite{Chen2019}.
In order to state  these results, we need some symbols for given graphs.

Let $S_{n,h}=K_h\nabla \overline{K}_{n-h}$.
Furthermore, $S^+_{n,h}=K_h \nabla (K_2\cup \overline{K}_{n-h-2})$.
Let $F_{n,k}=K_{k-1}\nabla ( (p K_2)\cup K_s)$, where $n-(k-1)=2p+s$ and $0\leq s<2$.
In addition, for $k\geq2$ and $d_1\geq\cdots \geq d_k\geq1$,  define
\begin{eqnarray*}
   f(k,d_1,\dots,d_k)=\frac{k^2(\sum_{i=1}^k d_i+k-2)^2(\sum_{i=1}^k 2d_i+5k-4)^{4k-2}+2(k-2)(\sum_{i=1}^kd_i)}{k-2}.
\end{eqnarray*}

  \begin{Theorem}\label{A}\cite{Chen2019}
 Let  $F=\cup_{i=1}^k P_{a_i}$  be a linear forest with   $k\geq2$ and
$a_1\geq \cdots \geq a_k\geq2$. Denote $ h=\sum \limits_{i=1}^k \lfloor\frac{a_i}{2}\rfloor-1$ and suppose that $G$ is  an $F$-free  graph  of sufficiently large order $n$.\\
(i) If there  exists an even $a_i$, then $\rho(G)\leq \rho(S_{n,h})$ with equality if and only if $G = S_{n,h}$;\\
(ii) If all $a_i$ are odd  and there exists at least one $a_i>3$,   then $\rho(G)\leq \rho(S^+_{n,h})$ with equality if and only if $G = S^+_{n,h}$.\\
(iii) If all $a_i$ are 3, i.e., $F= k P_3$, then $\rho(G)\leq \rho(F_{n,k}) $ with equality if and
only if   $G = F_{n,k}$.
\end{Theorem}

\begin{Theorem}\label{bipartite}\cite{Chen2019}
Let  $G$ be a $k P_3$-free bipartite graph of order $n\geq 11k-4$ with $k\geq2 $. Then $$\rho(G)\leq \sqrt{(k-1)(n-k+1)}$$ with equality if and
only if   $G= K_{k-1,n-k+1}$.
\end{Theorem}

\begin{Corollary}\label{least eigenvalue}\cite{Chen2019}
Let $G$ be a $kP_3$-free  graph of order $n\geq 11k-4$ with $k\geq2 $. Then $$\rho_n(G)\geq -\sqrt{(k-1)(n-k+1)}$$
 with equality if and only if   $G= K_{k-1,n-k+1}$.
\end{Corollary}

In Theorem~\ref{A}, the  extremal graph for $kP_3$ varies form other linear forests.  Note that $kP_3$ is also a  star forest $kS_2$.
Motivated by Problem~\ref{P1}, Theorems~\ref{thm-edge}, \ref{A} and \ref{bipartite},  we  determine the  maximum spectral radius and all  extremal graphs for all (bipartite) graphs of  order $n$ without  a star forest.   As a corollary, we determine the minimum  least eigenvalue of $A(G)$  and    all  extremal graphs for graphs  of  order $n$ without  a   star forest, extending Corollary~\ref{least eigenvalue} for large $n$. The main results of this paper are stated as follows.

\begin{Theorem}\label{thm1}
Let $F=\cup_{i=1}^k S_{d_i}$ be a star forest with $k\geq2$ and $d_1\geq\cdots \geq d_k\geq1$.
If $G$ be an $F$-free graph of  order  $n\geq \frac{(\sum_{i=1}^k 2d_i+5k-8)^4(\sum_{i=1}^k d_i+k-2)^4}{k-2}$,  then\\
 $$\rho(G)\leq \frac{k+d_k-3+\sqrt{(k-d_k-1)^2+4(k-1)(n-k+1)}}{2}$$ with equality if and only if $G=K_{k-1}\nabla H$, where $H$ is a $(d_k-1)$-regular graph of order $n-k+1$.
 In particular, if $d_k=2$, then $$\rho(G)\leq \rho(F_{n,k})$$ with equality if and only if $G=F_{n,k}$.
 \end{Theorem}


\vspace{3mm}
{\bf Remark 1.}  The extremal graph in Theorem~\ref{thm1}  only depends on the number of the components of $F$ and the minimum order of the stars in  $F$.

\begin{Theorem}\label{thm2}
Let $F=\cup_{i=1}^k S_{d_i}$ be a star forest with $k\geq2$ and $d_1\geq\cdots \geq d_k\geq1$. If $G$ is an $F$-free bipartite graph of order $n\geq \frac{f^2(k,d_1,\dots, d_k)}{4k-8}$, then
$$\rho(G)\leq \sqrt{(k-1)(n-k+1)}$$ with equality if and only if $G=K_{k-1,n-k+1}$.
\end{Theorem}

\begin{Corollary}\label{Cor3}
Let $F=\cup_{i=1}^k S_{d_i}$ be a star forest with $k\geq2$ and $d_1\geq\cdots \geq d_k\geq1$.  If $G$ is an $F$-free  graph of order $n\geq \frac{f^2(k,d_1,\dots, d_k)}{4k-8}$, then
$$\rho_n(G)\geq -\sqrt{(k-1)(n-k+1)}$$ with equality if and only if $G=K_{k-1,n-k+1}$.
\end{Corollary}

\vspace{3mm}
{\bf Remark 2.} For sufficiently large $n$, the extremal graphs in Theorem~\ref{thm2} and Corollary~\ref{Cor3} only depend on the number of the components of $F$.

\section{Preliminary}


%
%
%

  We first give a very rough estimation on the number of edges for a graph of order
$n\geq  \sum_{i=1}^k d_i+k$ without a star forest.
\begin{Lemma}\label{e1}
Let $F=\cup_{i=1}^k S_{d_i}$ be a star forest with $k\geq2$ and $d_1\geq\cdots \geq d_k\geq1$.
If $G$ is an $F$-free graph of order $n\geq  \sum_{i=1}^k d_i+k$, then
$$e(G)\leq\bigg(\sum_{i=1}^k d_i+2k-3\bigg)n-(k-1)\bigg(\sum_{i=1}^k d_i+k-1\bigg).$$
\end{Lemma}

\begin{Proof}
Let $C=\{v\in V(G):d(v)\geq \sum_{i=1}^k d_i+k-1\}$.  Since $G$ is $F$-free, $|C|\leq k-1$, otherwise we can embed an $F$ in $G$ by the definition of $C$.
Hence
\begin{eqnarray*}
  e(G) &=& \sum_{v\in C}d(v)+ \sum_{v\in V(G)\backslash C}d(v)\\
   &\leq&(n-1)|C|+(n-|C|)\bigg(\sum_{i=1}^k d_i+k-2\bigg) \\
   &=& \bigg(n- \sum_{i=1}^k d_i-k+1\bigg)|C|+\bigg(\sum_{i=1}^k d_i+k-2\bigg)n\\
   &\leq&(k-1)\bigg(n- \sum_{i=1}^k d_i-k+1\bigg)+\bigg(\sum_{i=1}^k d_i+k-2\bigg)n\\
   &=&\bigg(\sum_{i=1}^k d_i+2k-3\bigg)n-(k-1)\bigg(\sum_{i=1}^k d_i+k-1\bigg)\\
\end{eqnarray*}
\end{Proof}

\begin{Lemma}\label{spec1}
Let $F=\cup_{i=1}^k S_{d_i}$ be a star forest with $k\geq2$ and $d_1\geq\cdots \geq d_k\geq1$.
Let $G$ be an $F$-free connected bipartite graph of  order  $n\geq\frac{d_1^2}{k-1}+k-1$ with the maximum spectral radius  $\rho(G)$  and $\mathbf x=(x_u)_{u\in V(G)}$ be a positive eigenvector of $\rho(G)$ such that $\max\{x_u:u\in V(G)\}=1$. Then
 $x_u\geq \frac{1}{\rho (G)}$ for all $u\in V(G)$.
\end{Lemma}

\begin{Proof}
Set for short $\rho=\rho (G)$.  Choose a vertex $w\in V(G)$ such that $x_w=1$.
 Since $K_{k-1,n-k+1}$ is $F$-free, we have
$$\rho\geq \rho(K_{k-1,n-k+1})=\sqrt{(k-1)(n-k+1)}.$$

If $u=w$, then $x_u=1\geq \frac{1}{\rho }$. So  we  next suppose that $u\neq w$. We consider the following two cases.

\vspace{2mm}
{\bf Case~1.} $u$ is adjacent to $w$.
By eigenequation of $A(G)$ on $u$, $$\rho x_u=\sum_{uv\in E(G)}x_v\geq x_w=1,$$
which implies that $$x_u\geq \frac{1}{\rho }.$$

\vspace{2mm}
{\bf Case~2.} $u$ is not adjacent to $w$.
Let $G_1$ be a graph obtained from $G$ by deleting all edges incident with $u$ and adding an edge $uw$.
Note that $uw$ is a pendant edge in $G_1$.

{\bf Claim.} $G_1$ is also $F$-free.

Suppose that $G_1$ contains an $F$ as a subgraph. Since $G$ is $F$-free and $G_1$ contains an $F$ as a subgraph, we have $uw\in E(F)$. Since $uw$ is a pendant edge in $G_1$, $w$ is a center of $F$ with $d_F(w)=d_j$, where $1\leq j\leq k$. Let $G_2$ be the subgraph of $G_1$ by deleting $w$ and all its neighbors  in $F$. Note that $G_2$ is also a subgraph of $G$. Since $G_1$ contains an $F$ as a subgraph, $G_2$ contains $\cup_{ i\neq j}S_{d_i}$ as a subgraph. By eigenequation of $G$ on $w$,
$$d(w)\geq\sum_{vw\in E(G)}x_v=\rho x_w=\rho\geq \sqrt{(k-1)(n-k+1)}\geq d_1\geq d_j.$$
This implies that $G$ contains an $F$ as a subgraph, a contradiction.
\vspace{1mm}

By Claim, $G_1$ is $F$-free. Then
\begin{eqnarray*}
  0 &\geq& \rho(G_1)-\rho\geq \frac{\mathbf x^{\mathrm{T}}A(G_1)\mathbf x}{\mathbf x^{\mathrm{T}}\mathbf x}- \frac{ {\bf x^T} A(G)\bf x}{\bf {x^T}\bf x}\\
   &=& \frac{2}{\mathbf x^{\mathrm{T}}\mathbf x} \Big(x_ux_w-x_u\sum_{uv\in E(G)}x_v\Big)\\
 &=&  \frac{2x_u}{\mathbf x^{\mathrm{T}}\mathbf x} \Big(1-\rho x_u \Big),
\end{eqnarray*}
which implies that $$x_u\geq \frac{1}{\rho }.$$
This completes the proof.
\end{Proof}

\begin{Lemma}\label{spec2}
Let $d \geq 1$, $k\geq1$,  $n\geq \frac{(d-1)^2+(k-1)^2}{k-1}$, and  $H$ be a graph of order $n -k+1$.  If $G= K_{k-1} \nabla H$ and
$\Delta(H) \leq d-1$, then
$$\rho(G)\leq \frac{k+d-3+\sqrt{(k-d-1)^2+4(k-1)(n-k+1)}}{2}$$
with equality if and only if $H$ is a $(d-1)$-regular graph.
\end{Lemma}
\begin{Proof}
If $d=1$, then $G=K_{k-1}\nabla \overline{K}_{n-k+1}$. it is easy to calculate that
$$\rho (K_{k-1}\nabla \overline{K}_{n-k+1})=\frac{k-2+\sqrt{(k-2)^2+4(k-1)(n-k+1)}}{2}.$$ Next suppose that $d\geq 2$. Let $u_1,u_2,\cdots,u_{k-1}$ be the vertex of $G$ corresponding to $K_{k-1}$ in the representation $G := K_{k-1} \nabla H$. Set for short $\rho=\rho(G)$ and let ${\bf x}=(x_v)_{v\in E(G)}$  be a positive
eigenvector of $\rho$. By symmetry,  $x_{u_1}=\cdots=x_{u_{k-1}}$.  Choose a vertex $v\in V (H)$ such that
$$x_v = \max_{ w\in V (H)} x_w.$$
By eigenequation of $A(G)$ on $u_1$ and $v$, we have

\begin{equation}\label{6}
\begin{aligned}
 \rho x_{u_1} &= (k-2)x_{u_1}+\sum_{uu_1\in V(H)}x_u\leq(k-2)x_{u_1}+(n-k+1)x_v
     \end{aligned}
 \end{equation}
 \begin{equation}\label{7}
\begin{aligned}
  \rho x_{v}  &\leq& (k-1)x_{u_1}+\sum_{uv\in E(H)}x_u \leq(k-1)x_{u_1}+(d-1)x_v,
    \end{aligned}
 \end{equation}
which implies that
\begin{eqnarray*}
   (\rho-k+2)x_{u_1}&\leq& (n-k+1)x_v \\
 (\rho-d+1)x_v &\leq& (k-1)x_{u_1}.
\end{eqnarray*}
Since$$\rho> \rho(K_{k-1})= k-2,$$ and $$\rho>\rho(K_{k-1,n-k+1})=\sqrt{(k-1)(n-k+1)}\geq d-1,$$  we have
   $$ \rho^2-(k+d-3)\rho+ (k-2)(d-1)-(k-1)(n-k+1)\leq0.$$
Hence $$\rho\leq \frac{k+d-3+\sqrt{(k-d-1)^2+4(k-1)(n-k+1)}}{2}.$$
If equality holds, then all equalities in (\ref{6}) and (\ref{7}) hold. So  $d(v)=k+d-2$ and $x_u = x_v$ for any vertex
$u \in V (H)$. Since for any $u\in V(H)$,
\begin{eqnarray*}
 \rho x_u&=& (k-1)x_{u_1}+ \sum_{uz\in E(H)}x_z \leq(k-1)x_{u_1}+ (d-1)x_v=\rho x_v,
\end{eqnarray*}
we have $d(u)=d(v)=d+k-2$. So $H$ is $(d-1)$-regular.
\end{Proof}

\section{Proof of Theorem~\ref{thm1}}
Before proving Theorem~\ref{thm1}, we first prove the following  important result for connected graphs without a star forest.
\begin{Theorem}\label{c-thm1}
Let $F=\cup_{i=1}^k S_{d_i}$ be a star forest with $k\geq2$ and $d_1\geq\cdots \geq d_k\geq1$.
 If $G$ be an $F$-free connected graph of order $n\geq(\sum_{i=1}^k 2d_i+5k-7)^2(\sum_{i=1}^k d_i+k-2)^2$, then
 $$\rho(G)\leq \frac{k+d_k-3+\sqrt{(k-d_k-1)^2+4(k-1)(n-k+1)}}{2}$$ with equality if and only if $G=K_{k-1}\nabla H$, where $H$ is a $(d_k-1)$-regular graph of order $n-k+1$.
 In particular, if $d_k=2$, then $$\rho(G)\leq \rho(F_{n,k})$$ with equality if and only if $G=F_{n,k}$.
\end{Theorem}
\begin{Proof}
Let $G$ be an $F$-free connected  graph  of   order  $n$ with the maximum spectral radius.
Set for short $V=V(G)$, $E=E(G)$, $A=A(G)$, and $\rho=\rho(G)$. Let $\mathbf x=(x_v)_{v\in V(G)}$ be a positive eigenvector of $\rho $ such that $$x_w=\max\{x_u:u\in V(G)\}=1.$$

Since $K_{k-1,n-k+1}$ is $F$-free, we have
$$\rho\geq \rho(K_{k-1,n-k+1})=\sqrt{(k-1)(n-k+1)}. $$
Let $L=\{v\in V: x_v> \epsilon\} $ and $S=\{v\in V: x_v\leq \epsilon\} $, where $\epsilon = \frac{1}{\sum_{i=1}^k 2d_i+5k-7}$.

\vspace{2mm}
{\bf Claim.} $|L|=k-1$.

If $|L|\neq k-1$, then $|L|\geq k$ or $|L|\leq k-2$.

First suppose that $|L|\geq k$.
By eigenequation of $A$ on any vertex $u\in L$, we have
$$\sum_{i=1}^kd_i+k-2\leq \frac{\sqrt{(k-1)(n-k+1)}}{\sum_{i=1}^k 2d_i+5k-7}=\sqrt{(k-1)(n-k+1)}\epsilon<\rho x_u=\sum_{uv\in E} x_v\leq d(u),$$
where the first inequality holds because  $n\geq(\sum_{i=1}^k 2d_i+5k-7)^2(\sum_{i=1}^k d_i+k-2)^2$. Hence
$$ d(u)\geq\sum_{i=1}^kd_i+k-1.$$
Then we can embed an $F$ with all centers in $L$ in $G$, a contradiction.

Next suppose that $|L|\leq k-2$.
Then $$e(L)\leq \binom{|L|}{2}\leq \frac{1}{2}(k-2)(k-3)$$
and $$e(L,S)\leq (k-2)(n-k+2).$$ In addition,
by Lemma~\ref{e1}, $$e(S)\leq e(G)\leq\bigg(\sum_{i=1}^k d_i+2k-3\bigg)n.$$
By eigenequation of $A^2$ on  $w$, we have
\begin{eqnarray*}
  (k-1)(n-k+1) &\leq&\rho^2=\rho^2x_w=\sum\limits_{vw\in E}\sum\limits_{uv\in E}x_u\leq\sum\limits_{uv\in E}(x_u+x_v)\\
   &=& \sum\limits_{uv\in E(L,S)}(x_u+x_v)+\sum\limits_{uv\in E(S)}(x_u+x_v) +\sum\limits_{uv\in E(L)}(x_u+x_v)\\
   &\leq& \sum\limits_{uv\in E(L,S)}(x_u+x_v)+ 2\epsilon e(S)+2e(L)\\
   &\leq& \sum\limits_{uv\in E(L,S)}(x_u+x_v)+2\epsilon\bigg(\sum_{i=1}^k d_i+2k-3\bigg)n+(k-2)(k-3)
\end{eqnarray*}
Hence
\begin{eqnarray*}
  \sum\limits_{uv\in E(L,S)}(x_u+x_v)&\geq& (k-1)(n-k+1)-2\epsilon\bigg(\sum_{i=1}^k d_i+2k-3\bigg)n-(k-2)(k-3).
\end{eqnarray*}
On the other hand, by the definition of $L$ and $S$, we have
$$\sum\limits_{uv\in E(L,S)}(x_u+x_v)\leq(1+\epsilon)e(L,S)\leq (1+\epsilon)(k-2)(n-k+2).$$
Thus $$ (1+\epsilon)(k-2)(n-k+2)\geq (k-1)(n-k+1)-2\epsilon\bigg(\sum_{i=1}^k d_i+2k-3\bigg)n-(k-2)(k-3),$$
which implies that
$$\bigg(\bigg(\sum_{i=1}^k 2d_i+5k-8\bigg)\epsilon-1\bigg)n\geq\epsilon(k-2)^2-(k^2-3k+3).$$
Since $\epsilon = \frac{1}{\sum_{i=1}^k 2d_i+5k-7}$, we have
\begin{eqnarray*}
  n &\leq&(k^2-3k+3)\bigg(\sum_{i=1}^k 2d_i+5k-8\bigg)\bigg(\sum_{i=1}^k 2d_i+5k-7\bigg)-\\
  &&(k-2)^2\bigg(\sum_{i=1}^k 2d_i+5k-8\bigg) \\
   &\leq&  \bigg(\sum_{i=1}^k 2d_i+5k-7\bigg)^2\bigg(\sum_{i=1}^k d_i+k-2\bigg)^2,
\end{eqnarray*}
a contradiction.
This proves the Claim.

By Claim, $|L|=k-1$ and thus $|S|=n-k+1$.  Then  the subgraph $H$ induced by $S$ in $G$ is $S_{d_k}$-free. Otherwise, we can embed $F$ in $G$ with $k-1$ centers in $L$ and a center in $S$ as $d(u)\geq\sum_{i=1}^kd_i+k-1$ for any $u\in L$,  a contradiction.  Now  $\Delta(H)\leq d_k-1$. Note that the resulting graph  obtained  from $G$ by adding all edges in $L$ and all edges with one end in $L$ and the other in $S$ are also $F$-free and its spectral radius   increases strictly. By the extremality of $G$, we have $G=K_{k-1}\nabla H$. By Lemma~\ref{spec2} and the extremality of $G$, it follows that $H$ is a $(d_k-1)$-regular graph and
$$\rho=\frac{k+d_k-3+\sqrt{(k-d_k-1)^2+4(k-1)(n-k+1)}}{2}.$$
In particular, if $d_k=2$ then  $\Delta(H)\leq 1$, i.e.,  $H=pK_2\cup qK_1$, where $2p+q=n-k+1$. By the extremality of $G$,  $G=F_{n,k}$.
This completes the proof.
\end{Proof}

\vspace{3mm}

\noindent{\bf Proof of Theorems~\ref{thm1}.}
Let $G$ be an $F$-free  graph  of   order  $n$ with the maximum spectral radius.

 If $G$ is connected, then the result follows directly from Theorem~\ref{c-thm1}.
Next we suppose that $G$ is not connected. Since $K_{k-1,n-k+1}$ is $F$-free, we have
$$\rho(G)\geq \sqrt{(k-1)(n-k+1)}.$$
Let $G_1$ be a component of $G$  such that $\rho(G_1)=\rho(G)$ and $n_1=|V(G_1)|$.
Then
\begin{eqnarray*}
  n_1-1 &\geq&  \rho(G_1)=\rho(G)\geq \sqrt{(k-1)(n-k+1)}\geq \sqrt{(k-2)n}\\
  &\geq& \bigg(\sum_{i=1}^k 2d_i+5k-8 \bigg)^2 \bigg(\sum_{i=1}^k d_i+k-2 \bigg)^2,
\end{eqnarray*}
which implies that
\begin{eqnarray*}
  n_1 &\geq& \bigg(\sum_{i=1}^k 2d_i+5k-8\bigg)^2\bigg(\sum_{i=1}^k d_i+k-2\bigg)^2+1.
\end{eqnarray*}
By Theorem~\ref{c-thm1} again,
\begin{eqnarray*}
 \rho(G) =\rho(G_1) &\leq&  \frac{k+d_k-3+\sqrt{(k-d_k-1)^2+4(k-1)(n_1-k+1)}}{2} \\
   &<& \frac{k+d_k-3+\sqrt{(k-d_k-1)^2+4(k-1)(n-k+1)}}{2}.
\end{eqnarray*}
 In particular, if $d_k=2$ then it follows from By Theorem~\ref{c-thm1} again, $$\rho(G)=\rho(G_1)= \rho (F_{n_1,k})<\rho (F_{n,k}).$$
 Hence the result follows.
\QEDB

\vspace{3mm}
Note that extremal graph in Theorem~\ref{A} (iii) also holds for signless Laplacian special radius $q(G)$ \cite{Chen2020}. We conjecture the extremal graph in Theorem~\ref{thm1} also holds for signless Laplacian spectral radius $q(G)$.

\begin{Conjecture}
Let $F=\cup_{i=1}^k S_{d_i}$ be a star forest with $k\geq2$ and $d_1\geq\cdots \geq d_k\geq1$.
If $G$ be an $F$-free graph of  large order  $n$,  then
 $$q(G)\leq   \frac{n+2k+2d_k-6+\sqrt{(n+2k-2d_k-2)^2-8(k-1)(k-d_k-1)}}{2}$$ with equality if and only if $G=K_{k-1}\nabla H$, where $H$ is a $(d_k-1)$-regular graph of order $n-k+1$. In particular, if $d_k=2$, then $$q(G)\leq q(F_{n,k})$$ with equality if and only if $G=F_{n,k}$.
\end{Conjecture}
\section{Proofs of Theorem~\ref{thm2} and Corollary~\ref{Cor3}}

Before proving Theorem~\ref{thm2} and Corollary~\ref{Cor3}, we first prove the following  important result for bipartite connected graphs without a star forest.
\begin{Theorem}\label{c-thm3}
Let $F=\cup_{i=1}^k S_{d_i}$ be a star forest with $k\geq2$ and $d_1\geq\cdots \geq d_k\geq1$. If $G$ is an $F$-free connected bipartite graph of order $n\geq f(k,d_1,\dots,d_k)$, then
$$\rho(G)\leq \sqrt{(k-1)(n-k+1)}$$ with equality if and only if $G=K_{k-1,n-k+1}$.
\end{Theorem}
\begin{Proof}
Let $G$ be an $F$-free connected bipartite graph  of    order  $n$ with the maximum spectral radius.
Set for short $V=V(G)$, $E=E(G)$, $A=A(G)$, and $\rho=\rho(G)$. Let $\mathbf x=(x_v)_{v\in V(G)}$ be a positive eigenvector of $\rho $ such that $$x_w=\max\{x_u:u\in V(G)\}=1.$$

Since $K_{k-1,n-k+1}$ is $F$-free, we have
\begin{equation}\label{8}
\rho\geq \rho(K_{k-1,n-k+1})=\sqrt{(k-1)(n-k+1)}.
\end{equation}
Let $L=\{v\in V: x_v> \epsilon\} $ and $S=\{v\in V: x_v\leq \epsilon\} $, where $$\frac{\sum_{i=1}^kd_i+k-2}{\sqrt{(k-1)(n-k+1)}}\leq \epsilon \leq \frac{1}{k\sum_{i=1}^k (2d_i+5k-4)^{2k-1}}\bigg(1-\frac{\sum_{i=1}^k d_i}{n}\bigg).$$

{\bf Claim~1.} $|L|\leq k-1$.

Suppose that $|L|\geq k$.
By eigenequation of $A$ on any vertex $u\in L$, we have
$$\sum_{i=1}^kd_i+k-2\leq\sqrt{(k-1)(n-k+1)}\epsilon<\rho x_u=\sum_{uv\in E} x_v\leq d(u).$$
Hence $$d(u)\geq \sum_{i=1}^kd_i+k-1.$$
Then we can embed an $F$ in $G$ with all centers in $L$, a contradiction.
This proves Claim~1.

Since $|L|\leq k-1$, we have $$e(L)\leq \binom{|L|}{2}\leq \frac{1}{2}(k-1)(k-2)$$ and $$e(L,S)\leq (k-1)(n-k+1).$$ In addition,
by Lemma~\ref{e1}, $$e(S)\leq e(G)\leq\bigg(\sum_{i=1}^k d_i+2k-3\bigg)n.$$

 We next show that for any vertex  in $L$ has large degree.

{\bf Claim~2.}
Let $u\in L$ and $x_u=1-\delta$. Then
$$ d(u)\geq \bigg(1-\bigg(\sum_{i=1}^k 2d_i+5k-5\bigg)(\delta+\epsilon)\bigg)n.$$

Let $B_u=\{v\in V: uv\notin E\}$.
We first sum of eigenvector over all vertices of $G$.
\begin{eqnarray*}
  \rho \sum_{v\in V}x_v &=& \sum_{v\in V} \rho x_v=\sum_{v\in V} \sum_{vz\in E } x_z=\sum_{v\in V} d(v) x_v\leq \sum_{v\in L}  d(v)x_v + \sum_{v\in S} d(v)x_v \\
  &\leq& \sum_{v\in L}  d(v)+\epsilon \sum_{v\in S} d(v)=  2e(L)+e(L,S)+\epsilon(2e(S)+e(L,S))\\
   &=&  2e(L)+2\epsilon e(S)+(1+\epsilon)e(L,S),
\end{eqnarray*}
which implies that
\begin{equation}\label{4}
 \sum_{v\in V}x_v\leq \frac{ 2e(L)+2\epsilon e(S)+(1+\epsilon)e(L,S)}{\rho}.
\end{equation}
Next we sum of eigenvector over all vertices in $B_u$ by E.q. (\ref{4}) and Lemma~\ref{spec1}. Since
\begin{equation*}\label{5}
  \begin{aligned}
 \frac{1}{  \rho }|B_u|&\leq \sum_{v\in B_u}x_v\leq\sum_{v\in V(G)}x_v-\sum_{uv\in E(G)} x_v=\sum_{v\in V(G)}x_v-\rho x_u\\
   &\leq  \frac{ 2e(L)+2\epsilon e(S)+(1+\epsilon)e(L,S)}{\rho}-\rho x_u,
     \end{aligned}
\end{equation*}
 we have
\begin{eqnarray*}
  |B_u| &\leq&  2e(L)+2\epsilon e(S)+(1+\epsilon)e(L,S)-\rho^2x_u \\
   &\leq&  2e(L)+2\epsilon e(S)+(1+\epsilon)e(L,S)-(k-1)(n-k+1)(1-\delta) \\
   &\leq&(k-1)(k-2)+2\epsilon \bigg(\sum_{i=1}^k d_i+2k-3\bigg)n+(1+\epsilon)(k-1)(n-k+1)-\\
   &&(k-1)(n-k+1)(1-\delta)  \\
 &=&\bigg(2\epsilon \bigg(\sum_{i=1}^k d_i+2k-3\bigg)+(\delta+\epsilon)(k-1)\bigg)n+(k-1)(k-2)-(\delta+\epsilon)(k-1)^2\\
   &\leq&  \bigg(\sum_{i=1}^k 2d_i+4k-6+(k-1)+1\bigg)(\delta+ \epsilon)n\\
   &=& \bigg(\sum_{i=1}^k 2d_i+5k-6\bigg)(\delta+ \epsilon)n,
\end{eqnarray*}
where the last second inequality holds since $(k-1)(k-2)\leq \epsilon n<(\delta+\epsilon)n$ by the definition of $\epsilon$ and $n$.
Hence $$d(u)\geq n-1-\bigg(\sum_{i=1}^k 2d_i+5k-6\bigg)(\delta+ \epsilon)n\geq\bigg(1-\bigg(\sum_{i=1}^k 2d_i+5k-5\bigg)(\delta+\epsilon)\bigg)n.$$
This completes Claim~2.

\vspace{5mm}
{\bf Claim~3.}
 Let $ 1 \leq s<k-1 $. Suppose that  there is a set $X$ of $s$ vertices such that $X=\{v\in V:x_v\geq 1-\eta ~\text{and} ~d(v)\geq(1 - \eta)n\} $. Then there exists a vertex $u\in L\backslash X$ such that $$ x_{u}\geq 1 - \bigg(\sum_{i=1}^k2d_i+5k-5\bigg)^2(\eta + \epsilon)$$ and $$d(u)\geq\bigg(1 - \bigg(\sum_{i=1}^k2d_i+5k-5\bigg)^2(\eta + \epsilon)\bigg)n .$$

By eigenequation of $A^2$ on  $w$, we have
\begin{eqnarray*}
 \rho^2&=&\rho^2x_w=\sum\limits_{vw\in E}\sum\limits_{uv\in E}x_u\leq\sum\limits_{uv\in E}(x_u+x_v)\\
   &=& \sum\limits_{uv\in E(S)}(x_u+x_v) +\sum\limits_{uv\in E(L)}(x_u+x_v)\sum\limits_{uv\in E(L,S)}(x_u+x_v)\\
   &\leq&  2\epsilon e(S)+2e(L)+\sum\limits_{uv\in E(L,S)}(x_u+x_v)\\
  & \leq& 2\epsilon e(S)+2e(L)+\epsilon e(L,S)+\sum_{\substack{ uv\in E( L\backslash X,S)\\u \in L\backslash X }}x_u+\sum_{\substack{uv\in E(L\cap X,S)\\u\in L\cap X}}x_u,
\end{eqnarray*}
which implies that
\begin{eqnarray*}
 &&\sum_{\substack{uv\in E(L\backslash X,S)\\ u\in L\backslash X}}x_u\\
 &\geq&  \rho^2- 2\epsilon e(S)-2e(L)-\epsilon e(L,S)-\sum_{\substack{uv\in E( L\cap X,S)\\u\in L\cap X}}x_u\\
       &\geq&(k-1)(n-k+1)-2\epsilon\bigg (\sum_{i=1}^k d_i+2k-3\bigg)n-(k-1)(k-2)-\\
       &&\epsilon(k-1)(n-k+1)-sn\\
       &=&\bigg(k-1-s-\epsilon \bigg(\sum_{i=1}^k2d_i+5k-7\bigg)\bigg)n-(k-1)(2k-3)+\epsilon(k-1)^2\\
       &\geq&\bigg(k-1-s-\epsilon \bigg(\sum_{i=1}^k2d_i+5k-7\bigg)\bigg)n-\epsilon n\\
       &=&\bigg(k-1-s-\epsilon \bigg(\sum_{i=1}^k2d_i+5k-6\bigg)\bigg)n,
    \end{eqnarray*}
    where the last third inequality holds since $(k-1)(2k-3)\leq \epsilon n$ by the definition of $\epsilon$ and $n$.
In addition,
\begin{eqnarray*}
  e(L\backslash X,S)&=& e(L,S)-e(L\cap X, S) \\
   &\leq & (k-1)(n-k+1)-s(1-\eta)n+\binom{s}{2} \\
 &\leq&(k-1-s(1-\eta))n-\bigg((k-1)^2-\binom{k-2}{2} \bigg)\\
 &\leq&(k-1-s(1-\eta))n.
\end{eqnarray*}
Let $$g(s)=\frac{k-1-s-\epsilon \bigg(\sum_{i=1}^k2d_i+5k-6\bigg)}{k-1-s(1-\eta)}.$$  It is easy to see that $g(s)$ is decreasing with respect to $1\leq s\leq k-2$.
Then \begin{eqnarray*}
  \frac{\sum_{\substack{uv\in E(L\backslash X,S)\\u\in L\backslash X}}x_u}{e(L\backslash X,S)} &\geq& g(s)\geq g(k-2)
     =\frac{1-\epsilon \bigg(\sum_{i=1}^k2d_i+5k-6\bigg)}{1+(k-2)\eta}\\
     &\geq&1-\bigg(\sum_{i=1}^k2d_i+5k-6\bigg)(\eta+\epsilon).
\end{eqnarray*}
Hence there exists a vertex $u\in L\backslash X$ such that $$ x_u\geq1-\bigg(\sum_{i=1}^k2d_i+5k-6\bigg)(\eta+\epsilon).$$
By Claim~2,
\begin{eqnarray*}
d(u) &\geq& \bigg(1-\bigg(\sum_{i=1}^k 2d_i+5k-5\bigg)\bigg(\bigg(\sum_{i=1}^k2d_i+5k-6\bigg)(\eta+\epsilon)+\epsilon\bigg)\bigg)n\\
 &\geq &\bigg(1-\bigg(\sum_{i=1}^k2d_i+5k-5\bigg)^2(\eta+\epsilon)\bigg)n
\end{eqnarray*}
This completes Claim~3.
%
%

\vspace{5mm}
{\bf Claim~4.} $|L|=k-1$. Furthermore, for all $u\in L$, $$x_u\geq 1-\bigg (\sum_{i=1}^k 2d_i+5k-4\bigg)^{2k-1}\epsilon$$ and $$d(u)\geq\bigg(1-\bigg(\sum_{i=1}^k 2d_i+5k-4\bigg)^{2k-1}\epsilon\bigg)n.$$

Note that $w\in L$ and $x_{w}=1$. By Claim~2, $$d(w)\geq \bigg(1-\bigg(\sum_{i=1}^k 2d_i+5k-5\bigg)\epsilon\bigg)n. $$
 Applying Claim~5 iteratively for $k-2$ times,  we can find a set $X\subseteq L\backslash \{w\}$ of $k-2$ vertices such that for any $u\in X$,
  \begin{eqnarray*}
    x_u &\geq& 1-\bigg(\sum_{j=1}^{k-2}\bigg(\sum_{i=1}^k 2d_i+5k-5\bigg)^{2j}+\bigg(\sum_{i=1}^k 2d_i+5k-5\bigg)^{2k-2}\bigg(\sum_{i=1}^k 2d_i+5k-4\bigg)\bigg)\epsilon\\
     &\geq& 1- \bigg (\sum_{i=1}^k 2d_i+5k-4\bigg)^{2k-1}\epsilon
  \end{eqnarray*}
  and
  $$d(u)\geq\bigg(1-\bigg(\sum_{i=1}^k 2d_i+5k-4\bigg)^{2k-1}\epsilon\bigg)n.$$
Noting $|L|\leq k-1$, we have $L=X\cup \{w\}$. Hence $|L|=k-1$. This proves Claim~4.

\vspace{2mm}
 Let $T$ be the common neighborhood of $L$ and $R=S\backslash T$.  By Claim~4, $$|L|=k-1$$ and $$|T|\geq\bigg(1-k\bigg(\sum_{i=1}^k 2d_i+5k-4\bigg)^{2k-1}\epsilon\bigg)n\geq \sum_{i=1}^k d_i.$$  Since $G$ is bipartite, $L$ and $T$ are both independent sets of $G$.

\vspace{2mm}
{\bf Claim~5.}  $R$ is empty.

Suppose that  $R$ is not empty, i.e., there is a vertex  $v\in R$. Then $v$ has at most $d_k-1$ neighbors in $S$, otherwise we can embed an $F$ in $G$.
Let $H$ be a graph obtained from $G$ by removing all edges incident with $v$ and then connecting $v$ to each vertex in $L$. Clearly,
 $H$ is still $F$-free. By the definition of $R$, $v$ can be adjacent to at most $k-2$ vertices in $L$. Let $u\in L$ be the vertex not adjacent to $v$. Then By Claims~4 and 5, we have
\begin{eqnarray*}
  \rho(H)- \rho &\geq& \frac{{\bf x}^T A(H){\bf x}}{ {\bf x}^T{\bf x}}-\frac{{\bf x}^T A{\bf x}}{ {\bf x}^T{\bf x}}\\
 &\geq&\frac{2x_v}{{\bf x}^T{\bf x}}\bigg( x_u-\sum_{\substack{uz\in E\\ z\in S}} x_z\bigg)\\
   &\geq&  \frac{2x_v}{{\bf x}^T{\bf x}}\bigg(1-\bigg (\sum_{i=1}^k 2d_i+5k-4\bigg)^{2k-1}\epsilon-(d_k-1)\epsilon\bigg)\\
   &=& \frac{2x_v}{{\bf x}^T{\bf x}}\bigg( 1-\bigg (\bigg (\sum_{i=1}^k 2d_i+5k-4\bigg)^{2k-1}+d_k-1\bigg )\epsilon\bigg )\\
   &>&0,
\end{eqnarray*}
Hence $\rho(H)>\rho$, a contradiction. This proves Claim~5.

%
\vspace{2mm}
By Claim~5, $S=T$. By the definition of $T$, we have $G=K_{k-1,n-k+1}$. This completes the proof.
\end{Proof}

\vspace{5mm}
\noindent {\bf Proof of Theorem~\ref{thm2}.} Let $G$ be an $F$-free  bipartite graph  of   order  $n$ with the maximum spectral radius.

 If $G$ is connected, then the result follows directly from Theorem~\ref{c-thm3}.
Next we suppose that $G$ is not connected. Since $K_{k-1,n-k+1}$ is $F$-free,
$$\rho(G)\geq\sqrt{ (k-1)(n-k+1)}.$$ Let $G_1$ be a component of $G$  such that $\rho(G_1)=\rho(G)$ and $n_1=|V(G_1)|$.
Note that $G$ is triangle-free. By Wilf theorem \cite[Theorem~2]{Wilf1986}, we have
\begin{eqnarray*}
\frac{  n_1^2}{4} &\geq&  \rho^2(G_1)=\rho(G)^2\geq (k-1)(n-k+1)\geq (k-2)n
  \geq\frac{f^2(k,d_1,\dots,d_k)}{4},
\end{eqnarray*}
which implies that
$$n_1\geq f(k,d_1,\dots,d_k) . $$
By Theorem~\ref{c-thm3} again,
\begin{eqnarray*}
 \rho(G_1) &\leq&  \sqrt{ (k-1)(n_1-k+1)}< \sqrt{ (k-1)(n-k+1)},
\end{eqnarray*}
 a contradiction. This completes the proof.
\QEDB

\vspace{3mm}
\noindent {\bf Proof of Corollary~\ref{Cor3}.} By a result of Favaron et al. \cite{FMS}, $\rho_n(G) \geq \rho_n(H)$  for some spanning bipartite subgraph $H$. Moreover, the equality holds if and only if $G= H$, which can be deduced by its original proof.
 By Theorem~\ref{thm2},
$$\rho(H)\leq \sqrt{(k-1)(n-k+1)}$$ with equality if and only if $H = K_{k-1,n-k+1}$.
Since  the spectrum of a bipartite graph is symmetric \cite{LP},  $$\rho_n(H)\geq-\sqrt{(k-1)(n-k+1)}$$ with equality if and only if $H= K_{k-1,n-k+1}$.  Thus
 we have $$ \rho_n(G)\geq -\sqrt{(k-1)(n-k+1)}$$
with equality if and only if $G= K_{k-1,n-k+1}$.
\QEDB

\end{document}